\documentclass[10pt]{amsart}
\usepackage{}
\usepackage{mathrsfs}
\usepackage{amsfonts} 
\usepackage{graphicx,latexsym,bm,amsmath,amssymb,verbatim,multicol,lscape}

\theoremstyle{definition}

\theoremstyle{remark}

\numberwithin{equation}{section}

\begin{document}
\title [Asymptotic behavior of LCM of consecutive reducible quadratic progression terms]
{Asymptotic behavior of the least common multiple of consecutive reducible quadratic
progression terms}
\author{Guoyou Qian}
\address{Mathematical College, Sichuan University, Chengdu 610064, P.R. China}
\email{qiangy1230@gmail.com, qiangy1230@163.com}
\author{Shaofang Hong$^*$}
\address{Mathematical College, Sichuan University, Chengdu 610064, P.R. China}
\email{sfhong@scu.edu.cn, s-f.hong@tom.com, hongsf02@yahoo.com }
\thanks{$^*$Hong is the corresponding author and was supported partially
by National Science Foundation of China Grant \#11371260. Qian was supported
partially by Postdoctoral Science Foundation of China Grant \#2013M530109}
\date{\today}%
\begin{abstract}  Let $l$ and $m$ be two integers with $l>m\ge 0$, and let
$f(x)$ be the product of two linear polynomials with integer coefficients.
In this paper, we show that $\log {\rm lcm}_{mn<i\le ln}\{f(i)\}=An+o(n)$,
where $A$ is a constant depending only on $l$, $m$ and $f$. 
\end{abstract}

\maketitle

\section{\bf Introduction}
The study of the least common multiple of consecutive positive
integers was first initiated by Chebyshev for a significant attempt
to prove prime number theorem. From  Chebyshev's well-known work
\cite{[Ch]}, one can easily get an equivalent of prime number
theorem which states that $\log {\rm lcm}(1, ..., n)\sim n$ as $n$
tends to infinity. Since then, this topic  received attentions of
many authors. Hanson \cite{[Ha]} and Nair \cite{[N]} got the upper
and lower bound of ${\rm lcm}_{1\le i\le n}\{i\}$, respectively.
Bateman, Kalb and Stenger \cite{[BKS]} gave an asymptotic formula of
$\log{\rm lcm}_{1\le i\le n}\{b+ai\}$ as $n$ tends to infinity, where
$a$ and $b$ are  coprime integers. Farhi \cite{[F]}, Hong and Feng
\cite{[HF]}, Hong and Yang \cite{[HY]}, Hong and Kominers \cite{[HK]},
Wu, Tan and Hong \cite{[WTH]} and Kane and Kominers \cite{[KK]}
obtained lower bounds of the least common multiple of the first $n$
arithmetic progression terms. Farhi and Kane \cite{[FK]} studied the
least common multiple of consecutive integers. Hong and Qian \cite{[HQ]}
obtained some results on the least common multiple of consecutive arithmetic
progression terms which was consequently extended in one direction by Qian, Tan
and Hong \cite{[QTH2]}. Hong, Qian and Tan \cite{[HQT]} got an asymptotic
formula of the least common multiple of a sequence of products of linear
polynomials. On the other hand, Farhi \cite{[F]} obtained a nontrivial
lower bound for the least common multiple of the quadratic sequence $
\{i^2+1\}_{i=1}^\infty $. Oon \cite{[O]} improved some of the Hong-Kominers
result and Farhi's lower bound. Hong, Luo, Qian and Wang \cite{[HLQW]}
extended Nair's and Oon's lower bound by giving a uniform lower bound.
Qian, Tan and Hong \cite{[QTH]} showed that for any given positive
integer $k$, we have $\log{\rm lcm}_{0\le i\le k}\{(n+i)^2+1\}\sim 2(k+1)
\log n$ as $n\rightarrow \infty $. Recently, Hong and Qian \cite{[HQ2]}
got some interesting results on the least common multiple of consecutive
quadratic progression terms.

Qian and Hong \cite{[QH]} investigated the asymptotic
behavior of the least common multiple of any consecutive arithmetic
progression terms. Let $l$ and $m$ be integers
with $l>m\ge 0$ and let $a\ge 1$ and $b$ be integers such that
$a+b\ge 1$ and $\gcd(a, b)=1$. It is proved in \cite{[QH]} that
$$\log {\rm lcm}_{mn<i\le ln}\{ ai+b\}
\sim \frac{an}{\varphi(a)}\sum_{r=1\atop \gcd(r,a)=1}^{a}B_r
$$
as $n\rightarrow \infty$, where
\begin{align}\label{eq: 1.1}
B_r:={\left\{
\begin{array}{rl}
\frac{l}{r}, &\text{if} \ l\ge \frac{(a+r)m}{r},\\
\sum_{i=0}^{\mathcal{K}-1}\frac{l-m}{r+ai}+\frac{l}{r+a\mathcal{K}},
&\text{if} \ l<\frac{(a+r)m}{r}
\end{array}
\right.}
\end{align}
with
$\mathcal{K}:=\big\lfloor\frac{al-(l-m)r}{a(l-m)}\big\rfloor$
and $\lfloor x\rfloor$ being the largest integer no more than $x$.

In this paper, we mainly concentrate on the asymptotic behavior of the least
common multiple of consecutive reducible quadratic progression terms. There are
two cases about the reducible quadratic progressions. The first case
is $f(x)=(ax+b)^2$ with $a\ge 1$ and $b$ being integers such that $a+b\ge 1$
and $\gcd(a, b)=1$. This case is easy to answer. Actually, by the main result
of \cite{[QH]}, we can derive immediately that
$$\log {\rm lcm}_{mn<i\le ln}\{ (ai+b)^2\}
\sim \frac{2an}{\varphi(a)}\sum_{r=1\atop \gcd(r,a)=1}^{a}B_r
$$
as $n\rightarrow \infty$, where $B_r$ is defined as in (1.1).

Our main goal in the present paper is to treat with the second case
that $f(x)=(a_1x+b_1)(a_2x+b_2)$ with $a_i, b_i\in \mathbb{N}^*$
and $\gcd(a_i, b_i)=1$ for $1\le i\le 2$ and $a_1b_2\ne a_2b_1$.
Let $\mathbb{N}$ be the set of nonnegative integers and
$\mathbb{N}^*:=\mathbb{N}\setminus \{0\}$. For any two positive integers $a$ and $b$,
let $\langle b\rangle_{a}$ denote the smallest positive integer congruent to $b$
modulo $a$ between $1$ and $a$. For any integer $t$, we define $S_t$ by
$S_t:=\{i\in \mathbb{N} : 0\le i\le t\}$. Clearly, $S_t$ is empty if $t$
is negative and so we can define $\sum_{i\in S_t}g(i):=0$ for any arithmetic
function $g$ if $t<0$. We define the following three 4-variable arithmetic functions:
\begin{align}\label{eq: 1.2}
g_r(x,y,z,w):=\Big\lfloor \frac{xyl+ym\langle zr\rangle_{x}-xl\langle wr\rangle_{y}}
{xy(l-m)}\Big\rfloor,
\end{align}
\begin{align}\label{eq: 1.3}
 h_r(x,y,z,w):=\Big\lfloor \frac{xm\langle wr\rangle_{y}-yl
\langle zr\rangle_{x}}{xy(l-m)}\Big\rfloor
\end{align}
and
\begin{align}\label{eq: 1.4}
\lambda_r(x,y,z,w):=&\sum_{i\in S_{g_r(x,y,z,w)}}\frac{xl}{\langle zr\rangle_{x}+xi}-
\sum_{i\in S_{g_r(x,y,z,w)-1}}\frac{ym}{\langle wr\rangle_{y}+yi}\\
\nonumber&+\sum_{i\in S_{h_r(x,y,z,w)}}
\Big(\frac{yl}{\langle wr\rangle_{y}+yi}- \frac{xm}{\langle zr\rangle_{x}+xi}
\Big).
\end{align}
We can now state the main result of this paper.

\noindent{\bf Theorem 1.1.} {\it Let $l$ and $m$ be fixed integers
with $l>m\ge 0$.  Let $f(x)=(a_1x+b_1)(a_2x+b_2)$,
where $a_i, b_i\in \mathbb{N}^*$
and $\gcd(a_i, b_i)=1$ for $1\le i\le 2$ and $a_1b_2\ne a_2b_1$.
Then
$$\log {\rm lcm}_{mn<i\le ln}\{ f(i)\}
=\frac{n}{\varphi(q)}\sum_{r=1\atop \gcd(r,q)=1}^{q}A_r+o(n),
$$
where $q={\rm lcm}(a_1, a_2)$ and}
\begin{align}\label{eq: 1.5}
A_r:=
{\left\{
\begin{array}{rl}
\lambda_r(a_1, a_2, b_1, b_2) \ {\it if} \ a_1\langle b_2r\rangle_{a_2}\ge
a_2\langle b_1r\rangle_{a_1};\\
\lambda_r(a_2, a_1, b_2, b_1) \ {\it if} \ a_1
\langle b_2r\rangle_{a_2}<a_2\langle b_1r\rangle_{a_1}.
\end{array}
\right.}
\end{align}

Note that Theorem 1.1 is still true if at least one of $b_1$ and $b_2$
is a negative integer.

The paper is organized as follows. In Section 2, we prove two lemmas which
are needed for the proof of Theorem 1.1. The final section will devote to
the proof of Theorem 1.1.

\section{\bf Two lemmas}

In this section, we show two lemmas which are needed in the proof of
Theorem 1.1. Throughout, we let
\begin{align}\label{eq: 4.1}
H_1:=\Big\lfloor\frac{a_1l-(l-m)\langle b_1r\rangle_{a_1}}{a_1(l-m)}\Big\rfloor
\end{align}
and
\begin{align}\label{eq: 4.2}
H_2:=\Big\lfloor\frac{a_2l-(l-m)\langle b_2r\rangle_{a_2}}{a_2(l-m)}\Big\rfloor.
\end{align}
As usual, for any prime number $p$, we let $v_{p}$ be the normalized $p$-adic
valuation on the set of positive integers. Namely, one has $v_p(a)=s$ if
$p^{s}\parallel a$. We begin with the following result.

\noindent{\bf Lemma 2.1.} {\it Let $l, m, q$ and $f(x)$ be defined
as in Theorem 1.1. Then
$$\log {\rm lcm}_{mn<i\le ln}\{ f(i)\}
=\sum_{r'=1\atop \gcd(r',q)=1}^{q}\sum_{p\in \mathcal{P}_{r'}}\log p+O\big(\sqrt{n}\big),
$$
where
\begin{align}\label{eq: 4.3}
\mathcal {P}_{r'}:=\Big\{\text{prime}\ p: \ &p\equiv r'\pmod q \
\text{and} \ p\in \Big(0, (l-m)n\Big]\bigcup\\
\nonumber & \bigg(\bigcup_{j=1}^2\bigcup_{i=0}^{H_j}\Big(\frac{a_jmn}{\langle
b_jr\rangle_{a_j}+a_ji}, \frac{a_jln}
{\langle b_jr\rangle_{a_j}+a_ji}\Big]\bigg)\Big\}
\end{align}
with $r$ being the unique integer satisfying $rr'\equiv1\pmod q$ and $1\le r\le q$.}
\begin{proof}
For simplicity, we define $L_{m, l}^{(f)}(n):={\rm lcm}_{mn<i\le ln}\{f(i)\},$
and let $P_{m, l}^{(f)}(n)$ be the set of all the prime factors of
$L_{m,l}^{(f)}(n)$ not dividing ${\rm lcm}(a_1b_2-a_2b_1, q)$.

We claim that if $p\in P_{m, l}^{(f)}(n)$
and $p|f(i)$ for some integer $mn<i\le ln$, then $p$ divides exactly
one of $a_1i+b_1$ and $a_2i+b_2$. Otherwise, we have $p|(a_1i+b_1)$ and $p|(a_2i+b_2)$.
It implies that $p\mid \big(a_1(a_2i+b_2)-a_2(a_1i+b_1)\big)=a_1b_2-a_2b_1$,
which is impossible since $p\nmid {\rm lcm}(a_1b_2-a_2b_1, q)$.
The claim is proved. But the number of prime factors of
${\rm lcm}(a_1b_2-a_2b_1, q)$ is finite. So we have
\begin{align}\label{eq: 4.4}
&\log L_{m, l}^{(f)}(n)=\log \Big(\prod_{p\in P_{m, l}^{(f)}(n)}p^{v_p(L_{m, l}^{(f)}(n))}
\prod_{p\not\in P_{m, l}^{(f)}(n)}p^{v_p(L_{m, l}^{(f)}(n))}\Big)\\
\nonumber&=\sum_{p\in P_{m, l}^{(f)}(n)}v_p(L_{m, l}^{(f)}(n))\log p+O\big(\log\big(f(ln)\big)\big)\\
\nonumber &=\sum_{p\in P_{m, l}^{(f)}(n)}\log p+
\sum_{p\in P_{m, l}^{(f)}(n)\atop v_p(L_{m, l}^{(f)}(n))\ge 2
}\big(v_p\big(L_{m, l}^{(f)}(n)\big)-1\big)\log p+O\big(\log n\big).\quad\quad\
\end{align}
If  $p\in P_{m, l}^{(f)}(n)$ and $v_p(L_{m, l}^{(f)}(n))\ge 2$, then $p^2| f(i)$ for
some integer $i$ with $mn<i\le ln$.
Hence  by the claim we obtain that $p^2|(a_1i+b_1)$ or $p^2|(a_2i+b_2)$, which implies that
$$p\le M_n:=\max\{\sqrt{a_1ln+b_1)}, \sqrt{a_2ln+b_2}\}\ll \sqrt{n}.$$

On the other hand, since $p^{v_p(L_{m, l}^{(f)}(n))}\le f(ln)$, it follows that
$$v_p(L_{m, l}^{(f)}(n))\le \frac{\log f(ln)}{\log p}\ll \frac{\log n}{\log p}.$$
Hence we get by the prime number theorem that
$$\sum_{p\in P_{m, l}^{(f)}(n)\atop v_p(L_{m, l}^{(f)}(n))\ge 2
}\big(v_p(L_{m, l}^{(f)}(n))-1\big)\log p\ll \sum_{p\le M_n}\frac{\log n}
{\log p}\log p\ll \sum_{p\le M_n}\log n\ll \frac{\sqrt{n}}{\log \sqrt{n}}\log n\ll \sqrt{n}.$$
It then follows from (\ref{eq: 4.4}) that
\begin{align}\label{eq: 4.5}
\log L_{m, l}^{(f)}(n)=\sum_{p\in P_{m, l}^{(f)}(n)}\log p+O(\sqrt{n})+O(\log n)
=\sum_{p\in P_{m, l}^{(f)}(n)}\log p+O(\sqrt{n}).
\end{align}

First, we give a characterization on the primes in the set
$P_{m, l}^{(f)}(n)$. By $T(q)$ we denote the set of all
positive integers  no more than $q$ that are relatively prime to
$q$. Then by the definition of $P_{m, l}^{(f)}(n)$, we
know that each prime in $P_{m, l}^{(f)}(n)$ is relatively prime to
$q$. So each prime $p\in P_{m, l}^{(f)}(n)$ is congruent to $r'$
modulo $q$ for some $r'\in T(q)$. For convenience, we let
\begin{align}\label{eq: 4.6}
{\mathcal Q}_{r'}:=\{p\in P_{m, l}^{(f)}(n): p\equiv r'\pmod q\}.
\end{align}
Thus we derive from (\ref{eq: 4.5}) that
\begin{align}\label{eq: 4.7}
\log L_{m, l}^{(f)}(n)=\sum_{r'\in T(q)}\sum_{p\in P_{m,
l}^{(f)}(n)\atop p\equiv r'\pmod q}\log p+O(\sqrt{n})=\sum_{r'\in
T(q)}\sum_{p\in {\mathcal Q}_{r'}}\log p+O(\sqrt{n}).
\end{align}

For any given $r'\in T(q)$, there is exactly one $r\in T(q)$ such
that $rr'\equiv 1\pmod q$. Thus for any given prime $p\equiv r'\pmod q$,
we have $\langle b_jr\rangle_{a_j}p\equiv \langle
b_jr\rangle_{a_j}r'\equiv b_jrr'\equiv b_j\pmod{a_j}$ for each $1\le
j\le 2$. Since $\gcd(p, a_j)=1$ for $j=1, 2$, we can deduce
that all the terms divisible by $p$ in
the arithmetic progression $\{a_ji+b_j\}_{i=1}^{\infty}$ must be of
the form $(a_jk+\langle b_jr\rangle_{a_j})p$,
where $k\in \mathbb{N}$.  It follows that for each $1\le j\le 2$ and
any prime $p\in {\mathcal Q}_{r'}$, we have that $p|(a_ji+b_j)$ for some $mn< i\le
ln$ if and only if there is an integer $i_j\ge 0$ so that
$a_jmn+b_j<(a_ji_j+\langle b_jr\rangle_{a_j})p\le a_jln+b_j$.
Therefore, a prime $p$ congruent to $r'$ modulo $q$ is in $P_{m,
l}^{(f)}(n)$ if and only if $ p\nmid (a_1b_2-a_2b_1)$ and either
$$\frac{a_1mn+b_1}{\langle b_1r\rangle_{a_1}+a_1i_1}<p\le \frac{a_1ln+b_1}{\langle b_1r\rangle_{a_1}+a_1i_1}
$$
for some $i_1\in \mathbb{N}$, or
$$\frac{a_2mn+b_2}{\langle b_2r\rangle_{a_2}+a_2i_2}<p\le \frac{a_2ln+b_2}
{\langle b_2r\rangle_{a_2}+a_2i_2}$$
for some $i_2\in \mathbb{N}$.
Thus we have by (\ref{eq: 4.6}) that
\begin{align}\label{eq: 4.8}
{\mathcal Q}_{r'}=\bigcup_{j=1}^2\bigcup_{i=0}^{\infty}\Big\{ \text{prime}\ p\equiv r'\pmod q:
\frac{a_jmn+b_j}{\langle b_jr\rangle_{a_j}+a_ji}<p&\le
\frac{a_jln+b_j} {\langle b_jr\rangle_{a_j}+a_ji}\\
\nonumber&\ \mbox{and}\ p\nmid (a_1b_2-a_2b_1)\Big\}.
\end{align}
To prove Lemma 2.1, we have to treat with the union on the right-hand side of (\ref{eq: 4.8}).

Since $\gcd(p, a_j)=1$ for any prime $p\equiv r'\pmod q$, then by Lemma 3.6 of \cite{[HQ]},
there is exactly one term divisible by $p$ in any $p$ consecutive terms of
the arithmetic progression $\{ a_ji+b_j\}_{i=1}^{\infty}$ for each $1\le j\le 2$.
Therefore, for any prime $p$ with $p\le (l-m)n$ and $p\equiv r'\pmod q$, there is at
least one term divisible by $p$ in the set $\{(a_1i+b_1)(a_2i+b_2)\}_{i=mn+1}^{ln}$.
Hence we have
\begin{align}\label{eq: Q_3}
\big\{\text{prime}\ p\equiv r'\pmod q:\  p\le (l-m)n\ \mbox{and}\ p\nmid (a_1b_2-a_2b_1)\big\}
\subseteq {\mathcal Q}_{r'}.
\end{align}

By (\ref{eq: 4.1}) and (\ref{eq: 4.2}), for $j=1, 2$, we have that
$$\frac{a_jln+b_j}{\langle b_jr\rangle_{a_j}+a_j(H_j+1)}
<(l-m)n < \frac{a_jln+b_j}{\langle b_jr\rangle_{a_j}+a_jH_j}
$$
for any positive integer $n$ with
$$n>n_0:=\Big\lfloor\frac{b_j}{(a_j(H_j+1)
+\langle b_jr\rangle_{a_j})(l-m)-a_jl}\Big\rfloor.$$
It then follows that for $j=1, 2$ and all
integers $i$ with $i>H_j$, we have
$$\frac{a_jmn+b_j}{\langle b_jr\rangle_{a_j}+a_ji}<\frac{a_jln+b_j}
{\langle b_jr\rangle_{a_j}+a_ji}<(l-m)n$$
for any positive integer $n>n_0$. So we can deduce that
\begin{align*}
&\bigcup_{j=1}^2\bigcup_{i=H_j+1}^{\infty}\Big\{\text{prime}\ p\equiv r'\pmod q:
\frac{a_jmn+b_j}{\langle b_jr\rangle_{a_j}+a_ji}<p\le \frac{a_jln+b_j}{\langle b_jr\rangle_{a_j}+a_ji}\
\mbox{and}\\& p\nmid (a_1b_2-a_2b_1)\Big\}
\subseteq \Big\{\text{prime}\ p\equiv r'\pmod q: p\le (l-m)n\ \mbox{and}\ p\nmid (a_1b_2-a_2b_1)\Big\}
\end{align*}
for any positive integer $n>n_0$. It then follows from (\ref{eq: 4.8}) and (\ref{eq: Q_3}) that
\begin{align}\label{eq: 4.9}
{\mathcal Q}_{r'}&=\Big(\bigcup_{j=1}^2\bigcup_{i=0}^{H_j}\Big\{ \text{prime}\ p\equiv r'\pmod q:
\frac{a_jmn+b_j}{\langle b_jr\rangle_{a_j}+a_ji}<p\le \frac{a_jln+b_j}
{\langle b_jr\rangle_{a_j}+a_ji}\Big\}\\
\nonumber&\bigcup \{\text{prime}\ p\equiv r'\pmod q: p\le (l-m)n\}\Big)
\setminus\{\mbox{prime}\ p:\ p\nmid (a_1b_2-a_2b_1)\}
\end{align}
for any positive integer $n>n_0$.

Comparing (\ref{eq: 4.3}) with (\ref{eq: 4.9}) if $n>n_0$ and comparing (\ref{eq: 4.3})
with (\ref{eq: 4.8}) if $n\le n_0$, we know
that there are at most finitely many primes in the union set
$({\mathcal Q}_{r'}\setminus\mathcal{P}_{r'})\cup (\mathcal{P}_{r'}\setminus {\mathcal Q}_{r'})$
for any positive integer $n$.
Therefore
\begin{align}\label{eq: 4.10}
\sum_{p\in {\mathcal Q}_{r'}}\log p=\sum_{p\in \mathcal{P}_{r'}}\log p+O(\log n).
\end{align}
By (\ref{eq: 4.7}) and (\ref{eq: 4.10}), the desired result follows immediately.
This concludes the proof of Lemma 2.1.
\end{proof}

By Lemma 2.1, to estimate $\log {\rm lcm}_{mn<i\le ln}\{f(i)\}$, it suffices to estimate
$\sum_{p\in \mathcal{P}_{r'}}\log p$ for each integer $r'$ satisfying $1\le r'\le q$ and $\gcd(r', q)=1$,
which will be done in the following.

\noindent{\bf Lemma 2.2.} {\it Let $r'$ and $r$ be any given integers such that $1\le r', r\le q$ and
$rr'\equiv1\pmod{q}$. If $a_1\langle b_2r\rangle_{a_2}\ge a_2\langle b_1r\rangle_{a_1}$, then
$$\sum_{p\in \mathcal{P}_{r'}}\log p=\frac{n}{\varphi(q)} \lambda_r(a_1,a_2,b_1,b_2)+o(n),
$$
where $\mathcal{P}_{r'}$ and $\lambda_r(a_1, a_2, b_1,b_2)$ are defined as in (2.3) and (1.4), respectively.}

\begin{proof}
Since $a_1\langle b_2r\rangle_{a_2} \ge a_2\langle b_1r\rangle_{a_1}$,
we have
\begin{align}\label{eq: Q1}
\frac{a_1ln}{\langle b_1r\rangle_{a_1}+a_1i}\ge \frac{a_2ln}
{\langle b_2r\rangle_{a_2}+a_2i}\ \mbox{and}\ \frac{a_1mn}{\langle b_1r\rangle_{a_1}+a_1i}\ge \frac{a_2mn}
{\langle b_2r\rangle_{a_2}+a_2i}
\end{align}
for any integer $i\ge 0$. On the other hand, for any integer $i\ge 0$, we have
\begin{align}\label{eq: Q2}
\frac{a_1ln}{\langle b_1r\rangle_{a_1}+a_1(i+1)}< \frac{a_2ln}
{\langle b_2r\rangle_{a_2}+a_2i}\ \mbox{and}\ \frac{a_1mn}{\langle b_1r\rangle_{a_1}+a_1(i+1)}\le \frac{a_2mn}
{\langle b_2r\rangle_{a_2}+a_2i}
\end{align}
since $0\le a_1\langle b_2r\rangle_{a_2} -a_2\langle b_1r\rangle_{a_1}< a_1a_2$ and $l>m\ge 0$.
Let $K_1=g_r(a_1, a_2, b_1, b_2)$ and $K_2=h_r(a_1, a_2, b_1, b_2)$.
Then by (\ref{eq: 1.2}) and (\ref{eq: 1.3}), we get
\begin{align}\label{eq: 4.13}
K_1=\Big\lfloor \frac{a_1a_2l+a_2\langle b_1r\rangle_{a_1}m-a_1\langle b_2r\rangle_{a_2}l}
{a_1a_2(l-m)}\Big\rfloor
\end{align}
and
\begin{align}\label{eq: 4.14}
K_2=\Big\lfloor \frac{a_1\langle b_2r\rangle_{a_2}m-a_2
\langle b_1r\rangle_{a_1}l}{a_1a_2(l-m)}\Big\rfloor.
\end{align}
Thus by (1.4),  in order to show Lemma 2.2, we only need to prove that
\begin{align}\label{eq: 4.15}
\sum_{p\in \mathcal{P}_{r'}}\log p=& \frac{n}{\varphi(q)}\bigg(\sum_{i\in S_{K_1}}
\frac{a_1l}{\langle b_1r\rangle_{a_1}+a_1i}-
\sum_{i\in S_{K_1-1}}\frac{a_2m}{\langle b_2r\rangle_{a_2}+a_2i}+\\
\nonumber \ &\sum_{i\in S_{K_2}}
\Big(\frac{a_2l}{\langle b_2r\rangle_{a_2}+a_2i}- \frac{a_1m}{\langle b_1r\rangle_{a_1}+a_1i}
 \Big)\bigg)+o(n).
\end{align}

In the following we show that (\ref{eq: 4.15}) is true. For this purpose,
we need to analyze the following union
\begin{align}\label{eq: 4.16}
\mathcal{T}_r:=\Big(\bigcup_{j=1}^2\bigcup_{i=0}^{H_j}\Big(\frac{a_jmn}{\langle
b_jr\rangle_{a_j}+a_ji}, \frac{a_jln}
{\langle b_jr\rangle_{a_j}+a_ji}\Big]\Big)\bigcup\big(0, (l-m)n\big],
\end{align}
since  (\ref{eq: 4.3}) gives that
\begin{align}
 \label{eq: 4.26}
\mathcal{P}_{r'}=\{\text{prime}\ p\equiv r'\pmod q:
\ p\in  \mathcal{T}_r\}.
\end{align}

Evidently, we have
\begin{align*}
 \frac{a_1l-(l-m)\langle b_1r\rangle_{a_1}}{a_1(l-m)}
-\frac{a_2l-(l-m)\langle b_2r\rangle_{a_2}}{a_2(l-m)}
=\frac{a_1\langle b_2r\rangle_{a_2}-a_2\langle b_1r\rangle_{a_1}}{a_1a_2}
\end{align*}
and $$0\le  \frac{a_1\langle b_2r\rangle_{a_2}-a_2\langle b_1r\rangle_{a_1}}{a_1a_2}<1.$$
Thus by  (\ref{eq: 4.1}) and  (\ref{eq: 4.2}) we get that
\begin{align}\label{eq: 4.17}
H_1=H_2\ \mbox{or}\  H_2+1.
\end{align}
Moreover, for each $1\le j\le 2$, it follows from (\ref{eq: 4.1}) and  (\ref{eq: 4.2}) that
$$\frac{a_jm-(l-m)\langle b_jr\rangle_{a_j}}{a_j(l-m)}<H_j\le
\frac{a_jl-(l-m)\langle b_jr\rangle_{a_j}}{a_j(l-m)}.$$
Hence for each $1\le j\le 2$,
\begin{align}\label{eq: 4.18}
\frac{a_jmn}{\langle b_jr\rangle_{a_j}+a_jH_j}<(l-m)n\le \frac{a_jln}
{\langle b_jr\rangle_{a_j}+a_jH_j}.
\end{align}
By (\ref{eq: 4.13}), we have $K_1\ge 0$ and
\begin{align*}
K_1-1\le  \frac{a_1a_2m+a_2\langle b_1r\rangle_{a_1}m-a_1
\langle b_2r\rangle_{a_2}l}{a_1a_2(l-m)}<K_1.
\end{align*}
It then follows that
\begin{align}\label{eq: 4.19}
\frac{a_1ln}{\langle b_1r\rangle_{a_1}+a_1(i+1)}> \frac{a_2mn}{\langle b_2r\rangle_{a_2}+a_2i}
\end{align}
for any $i\ge K_1$ and
\begin{align}\label{eq: 4.20}
\frac{a_1ln}{\langle b_1r\rangle_{a_1}+a_1(i+1)}\le \frac{a_2mn}{\langle b_2r\rangle_{a_2}+a_2i}
\end{align}
for any $0\le i\le K_1-1$ if $K_1\ge 1$.

From (\ref{eq: 4.14}), we know that $K_2$ may be smaller than 0, and
$$
a_1a_2K_2(l-m)\le a_1\langle b_2r\rangle_{a_2}m-a_2
\langle b_1r\rangle_{a_1}l<a_1a_2(K_2+1)(l-m).
$$
Thus
\begin{align}\label{eq: 4.21}
\frac{a_2ln}{\langle b_2r\rangle_{a_2}+a_2i}
>\frac{a_1mn}{\langle b_1r\rangle_{a_1}+a_1i}.
\end{align}
for any $i\ge \max(0, K_2+1)$, and
\begin{align}\label{eq: 4.22}
\frac{a_2ln}{\langle b_2r\rangle_{a_2}+a_2i}
\le \frac{a_1mn}{\langle b_1r\rangle_{a_1}+a_1i}
\end{align}
for any  $0\le i\le K_2$ if $K_2\ge 0$.

For $j=1, 2$, if $H_j\ge 1$, then by (\ref{eq: 4.1}) and (\ref{eq: 4.2}) we infer that
\begin{align*}
\frac{a_jmn}{\langle b_jr\rangle_{a_j}+a_j(i-1)}-\frac{a_jln}{\langle
b_jr\rangle_{a_j}+a_ji}&=
\frac{a_j\big(a_jln-(l-m)\langle b_jr\rangle_{a_j}n-a_ji(l-m)n\big)}
{(\langle b_jr\rangle_{a_j}+a_j(i-1))(\langle b_jr\rangle_{a_j}+ai)}\\
&\ge \frac{a_j(a_jH_j(l-m)n-a_ji(l-m)n)}{(\langle b_jr\rangle_{a_j}+a_j(i-1))(\langle b_jr\rangle_{a_j}+a_ji)}\\
&= \frac{a_j^2(H_j-i)(l-m)n}{(\langle
b_jr\rangle_{a_j}+a_j(i-1))(\langle b_jr\rangle_{a_j}+a_ji)}\ge 0
\end{align*}
for any integer $i$ with $1\le
i\le H_j$, which means that
$$\frac{a_jmn}{\langle b_jr\rangle_{a_j}+a_j(i-1)}\ge \frac{a_jln}{\langle
b_jr\rangle_{a_j}+a_ji}$$
for any integer $i$ with $1\le i\le H_j$.
Hence for $j=1, 2$, the intersection
\begin{align}\label{eq: 4.23}
\Big(\frac{a_jmn}{\langle b_jr\rangle_{a_j}+a_ji_1},
\frac{a_jln}{\langle b_jr\rangle_{a_j}+a_ji_1}\Big] \bigcap \Big( \frac{a_jmn}{\langle b_jr\rangle_{a_j}+a_ji_2},
\frac{a_jln}{\langle b_jr\rangle_{a_j}+a_ji_2}\Big]
\end{align}
is empty for any $0\le i_1\ne i_2\le H_j$ if $H_j\ge 1$. Now we consider the following two cases.

{\sc Case 1.}  $K_1\ge K_2+1$.  First, it is easy to see from (\ref{eq: 4.2}), (\ref{eq: 4.13})
and (\ref{eq: 4.14}) that $K_1\le H_2$ and $K_2+1\le H_2$. For any integer $i\ge \max(0, K_2+1)$,
we have by (\ref{eq: Q1}) and (\ref{eq: 4.21}) that
\begin{align}\label{eq: 4.24}
\bigcup_{j=1}^2\Big( \frac{a_jmn}{\langle b_jr\rangle_{a_j}+a_ji},
\frac{a_jln}{\langle b_jr\rangle_{a_j}+a_ji}\Big]=\Big(\frac{a_2mn}{\langle b_2r\rangle_{a_2}+a_2i},
\frac{a_1ln}{\langle b_1r\rangle_{a_1}+a_1i}\Big].
\end{align}
It then follows from (\ref{eq: 4.17})-(\ref{eq: 4.19}) and (\ref{eq: 4.24}) that
\begin{align}\label{eq: 4.25}
&\bigg(\bigcup_{j=1}^2\bigcup_{i=K_1}^{H_j}\Big( \frac{a_jmn}{\langle b_jr\rangle_{a_j}+a_ji},
\frac{a_jln}{\langle b_jr\rangle_{a_j}+a_ji}\Big]\bigg)\bigcup \Big( 0, (l-m)n\Big]\\
\nonumber&=\bigg(\bigcup_{i=K_1}^{H_2}\Big( \frac{a_2mn}{\langle b_2r\rangle_{a_2}+a_2i},
\frac{a_1ln}{\langle b_1r\rangle_{a_1}+a_1i}\Big]\bigg) \bigcup \Big( 0, (l-m)n\Big]\\
\nonumber &\quad \ \bigcup\Big( \frac{a_1mn}{\langle b_1r\rangle_{a_1}+a_1H_1},
\frac{a_1ln}{\langle b_1r\rangle_{a_1}+a_1H_1}\Big]\\
\nonumber &= \Big( \frac{a_2mn}{\langle b_2r\rangle_{a_2}+a_2H_2},
\frac{a_1ln}{\langle b_1r\rangle_{a_1}+a_1K_1}\Big] \bigcup \Big( 0,
\frac{a_1ln}{\langle b_1r\rangle_{a_1}+a_1H_1}\Big]\\
\nonumber &=\Big( 0,
\frac{a_1ln}{\langle b_1r\rangle_{a_1}+a_1K_1}\Big].
\end{align}
Thus we can derive from  (\ref{eq: 4.16}) and (\ref{eq: 4.25}) that

\begin{align*}
\mathcal{T}_r=&\bigcup_{j=1}^2\bigg(\bigcup_{i\in S_{K_1-1}}\Big( \frac{a_jmn}{\langle b_jr\rangle_{a_j}+a_ji},
\frac{a_jln}{\langle b_jr\rangle_{a_j}+a_ji}\Big]\bigcup \\
\nonumber&\bigcup_{i=K_1}^{H_j}\Big( \frac{a_jmn}{\langle b_jr\rangle_{a_j}+a_ji},
\frac{a_jln}{\langle b_jr\rangle_{a_j}+a_ji}\Big]\bigg)
 \bigcup \Big( 0, (l-m)n\Big]\\
\nonumber =&\bigg(\bigcup_{j=1}^2\bigcup_{i\in S_{K_2}}\Big( \frac{a_jmn}{\langle b_jr\rangle_{a_j}+a_ji},
\frac{a_jln}{\langle b_jr\rangle_{a_j}+a_ji}\Big]\bigg)\bigcup\Big(0,
\frac{a_1ln}{\langle b_1r\rangle_{a_1}+a_1K_1}\Big]\\
\nonumber &\bigcup  \bigg(\bigcup_{i\in S_{K_1-1}\setminus
S_{K_2}}\bigcup_{j=1}^2\Big( \frac{a_jmn}{\langle b_jr\rangle_{a_j}+a_ji},
\frac{a_jln}{\langle b_jr\rangle_{a_j}+a_ji}\Big]\bigg).
\end{align*}
It then follows from (\ref{eq: 4.24}) that
\begin{align}\label{eq: 4.27}
\mathcal{T}_r=&\bigg(\bigcup_{j=1}^2\bigcup_{i\in S_{K_2}}
\Big( \frac{a_jmn}{\langle b_jr\rangle_{a_j}+a_ji},
\frac{a_jln}{\langle b_jr\rangle_{a_j}+a_ji}\Big]\bigg)\bigcup\\
\nonumber &\bigg(\bigcup_{i\in S_{K_1-1}\setminus S_{K_2}}
\Big(\frac{a_2mn}{\langle b_2r\rangle_{a_2}+a_2i},
\frac{a_1ln}{\langle b_1r\rangle_{a_1}+a_1i}\Big]\bigg)\bigcup\Big(0,
\frac{a_1ln}{\langle b_1r\rangle_{a_1}+a_1K_1}\Big].
\end{align}

Note that $S_{K_2}$ is empty if $K_2<0$, and $S_{K_1-1}\setminus S_{K_2}$ is empty
if $K_1=K_2+1$ or $K_1=0$. By (\ref{eq: 4.20}), we know that the following union
$$
\bigcup_{i\in S_{K_1-1}\setminus S_{K_2}}\Big(\frac{a_2mn}{\langle b_2r\rangle_{a_2}+a_2i},
\frac{a_1ln}{\langle b_1r\rangle_{a_1}+a_1i}\Big]
$$
is a disjoint union. But by (\ref{eq: 4.20}), (\ref{eq: 4.22}) and (\ref{eq: 4.23}), the union
$$
\bigcup_{j=1}^2\bigcup_{i\in S_{K_2}}\Big( \frac{a_jmn}{\langle b_jr\rangle_{a_j}+a_ji},
\frac{a_jln}{\langle b_jr\rangle_{a_j}+a_ji}\Big]
$$
is a disjoint union. Therefore by (\ref{eq: 4.20}), the union on the right-hand side
of (\ref{eq: 4.27}) is disjoint. Thus applying (\ref{eq: 4.26}), (\ref{eq: 4.27})
and prime number theorem for arithmetic progressions (see, for example \cite{[MV]}),
we obtain that
\begin{align*}
\sum_{p\in \mathcal{P}_{r'}}\log p&=\sum_{j=1}^2\sum_{i\in S_{K_2}}
\sum_{\frac{a_jmn}{\langle b_jr\rangle_{a_j}+a_ji}<
p\le \frac{a_jln}{\langle b_jr\rangle_{a_j}+a_ji}\atop p\equiv r'\pmod q}\log p
+\sum_{p\le \frac{a_1ln}{\langle b_1r\rangle_{a_1}+a_1K_1}\atop p\equiv r'\pmod q}\log p\\
&\ \ +\sum_{i\in S_{K_1-1}\setminus S_{K_2}}\sum_{
\frac{a_2mn}{\langle b_2r\rangle_{a_2}+a_2i}<p\le \frac{a_1ln}{\langle b_1r
\rangle_{a_1}+a_1i}\atop p\equiv r'\pmod q}\log p\\
&=\frac{n}{\varphi(q)}\bigg( \sum_{j=1}^2\sum_{i\in S_{K_2}}\Big( \frac{a_jl}{\langle b_jr\rangle_{a_j}+a_ji}-
\frac{a_jm}{\langle b_jr\rangle_{a_j}+a_ji} \Big)+ \frac{a_1l}{\langle b_1r\rangle_{a_1}+a_1K_1} \\
&\quad\quad\quad\quad+\sum_{i\in S_{K_1-1}\setminus S_{K_2}}\Big( \frac{a_1l}{\langle b_1r\rangle_{a_1}+a_1i}-
\frac{a_2m}{\langle b_2r\rangle_{a_2}+a_2i}\Big)\bigg)+o(n).
\end{align*}
Then (\ref{eq: 4.15}) follows immediately. So (\ref{eq: 4.15}) is proved for Case 1.

{\sc Case 2.}  $K_1\le K_2$.  Then by (\ref {eq: 4.13}) and (\ref {eq: 4.14}), we have $K_2\ge K_1\ge 0$.
If $K_2+1\le H_2$, applying (\ref{eq: 4.17})-(\ref{eq: 4.19}) and (\ref{eq: 4.24}),
one infers that

\begin{align} \label{eq: 4.30}
&\bigcup_{j=1}^2\bigcup_{i=K_2+1}^{H_j}\Big(\frac{a_jmn}{\langle b_jr\rangle_{a_j}+a_ji},
\frac{a_jln}{\langle b_jr\rangle_{a_j}+a_ji}\Big]\bigcup\Big( 0, (l-m)n\Big]\\
\nonumber &=\bigcup_{i=K_2+1}^{H_2}\Big(\frac{a_2mn}{\langle b_2r\rangle_{a_2}+a_2i},
\frac{a_1ln}{\langle b_1r\rangle_{a_1}+a_1i}\Big]\bigcup\Big( 0, (l-m)n\Big]\\
\nonumber &\quad\quad\bigcup \Big(\frac{a_1mn}{\langle b_1r\rangle_{a_1}+a_1H_1},
\frac{a_1ln}{\langle b_1r\rangle_{a_1}+a_1H_1}\Big]\\
\nonumber &=\Big(\frac{a_2mn}{\langle b_2r\rangle_{a_2}+a_2H_2},
\frac{a_1ln}{\langle b_1r\rangle_{a_1}+a_1(K_2+1)}\Big]
\bigcup\Big(0, \frac{a_1ln}{\langle b_1r\rangle_{a_1}+a_1H_1}\Big]\\
\nonumber &=\Big(0, \frac{a_1ln}{\langle b_1r\rangle_{a_1}+a_1(K_2+1)}\Big].
\end{align}
Hence by (\ref{eq: 4.16}) and (\ref{eq: 4.30}),
\begin{align}\label{eq: 4.31}
\mathcal{T}_r&=\bigg(\bigcup_{j=1}^2\bigcup_{i=0}^{K_2}
\Big(\frac{a_jmn}{\langle b_jr\rangle_{a_j}+a_ji},
\frac{a_jln}{\langle b_jr\rangle_{a_j}+a_ji}\Big]\bigg)\bigcup
\Big(0, \frac{a_1ln}{\langle b_1r\rangle_{a_1}+a_1(K_2+1)}\Big].
\end{align}
Moreover, we have

\begin{align}
\label{eq: 4.32}
&\bigcup_{j=1}^2\bigcup_{i=0}^{K_2}\Big(\frac{a_jmn}{\langle b_jr\rangle_{a_j}+a_ji},
\frac{a_jln}{\langle b_jr\rangle_{a_j}+a_ji}\Big]\\
\nonumber &=\bigg(\bigcup_{j=1}^2\bigcup_{i\in S_{K_1-1}}
\Big(\frac{a_jmn}{\langle b_jr\rangle_{a_j}+a_ji},
\frac{a_jln}{\langle b_jr\rangle_{a_j}+a_ji}\Big]\bigg)\bigcup
\Big(\frac{a_1mn}{\langle b_1r\rangle_{a_1}+a_1K_1},
\frac{a_1ln}{\langle b_1r\rangle_{a_1}+a_1K_1}\Big] \bigcup\\
\nonumber &\bigcup_{i\in S_{K_2-1\setminus S_{K_1-1}}}
\bigg( \Big(\frac{a_1mn}{\langle b_1r\rangle_{a_1}+a_1(i+1)},
\frac{a_1ln}{\langle b_1r\rangle_{a_1}+a_1(i+1)}\Big]\bigcup\\
\nonumber &  \Big(\frac{a_2mn}{\langle b_2r\rangle_{a_2}+a_2i},
\frac{a_2ln}{\langle b_2r\rangle_{a_2}+a_2i}\Big]\bigg)
\bigcup \Big(\frac{a_2mn}{\langle b_2r\rangle_{a_2}+a_2K_2},
\frac{a_2ln}{\langle b_2r\rangle_{a_2}+a_2K_2}\Big].
\end{align}
But since $K_2\ge 0$ and $K_2\ge K_1$, applying (\ref{eq: Q2}) and (\ref {eq: 4.19}) gives us that
\begin{align}\label{eq: 4.33}
\Big(0, \frac{a_1ln}{\langle b_1r\rangle_{a_1}+a_1(K_2+1)}\Big] \bigcup
\Big(\frac{a_2mn}{\langle b_2r\rangle_{a_2}+a_2K_2},
\frac{a_2ln}{\langle b_2r\rangle_{a_2}+a_2K_2}\Big]
=\Big(0, \frac{a_2ln}{\langle b_2r\rangle_{a_2}+a_2K_2}\Big].
\end{align}
Therefore by (\ref{eq: Q2}), (\ref{eq: 4.19}) and (\ref{eq: 4.31})-(\ref{eq: 4.33}), we have
\begin{align}\label{eq: 4.34}
\mathcal{T}_r&=\bigg(\bigcup_{j=1}^2\bigcup_{i\in S_{K_1-1}}\Big(\frac{a_jmn}{\langle b_jr\rangle_{a_j}+a_ji},
\frac{a_jln}{\langle b_jr\rangle_{a_j}+a_ji}\Big]\bigg)\bigcup \Big(\frac{a_1mn}{\langle b_1r\rangle_{a_1}+a_1K_1},
\frac{a_1ln}{\langle b_1r\rangle_{a_1}+a_1K_1}\Big]\\
\nonumber &\bigg(\bigcup_{i\in S_{K_2-1\setminus S_{K_1-1}}} \Big(\frac{a_1mn}{\langle b_1r\rangle_{a_1}+a_1(i+1)},
\frac{a_2ln}{\langle b_2r\rangle_{a_2}+a_2i}\Big]\bigg)\bigcup \Big(0, \frac{a_2ln}{\langle b_2r\rangle_{a_2}+a_2K_2}\Big].
\end{align}

Note that $S_{K_2-1}\setminus S_{K_1-1}$ is empty if $K_1=K_2$, and $S_{K_1-1}$ is empty if $K_1=0$.
By (\ref{eq: 4.20}), we have that
$$\frac{a_1ln}{\langle b_1r\rangle_{a_1}+a_1K_1}\le \frac{a_2mn}{\langle b_2r\rangle_{a_2}+a_2(K_1-1)}$$
if $K_1\ge 1$.
Using (\ref{eq: 4.22}), one has
$$
\frac{a_2ln}{\langle b_2r\rangle_{a_2}+a_2K_2}\le \frac{a_1mn}{\langle b_1r\rangle_{a_1}+a_1K_2}
\ {\rm and} \
\frac{a_2ln}{\langle b_2r\rangle_{a_2}+a_2K_1}\le \frac{a_1mn}{\langle b_1r\rangle_{a_1}+a_1K_1}.
$$
Then using (\ref{eq: 4.20}), (\ref{eq: 4.22}) and (\ref{eq: 4.23}), we derive that any two intervals in the
union on the right-hand side of (\ref{eq: 4.34}) are disjoint.
Hence we have by (\ref{eq: 4.26}) and (\ref{eq: 4.34}) that
\begin{align*}
&\sum_{p\in \mathcal{P}_{r'}}\log p=\sum_{j=1}^2\sum_{i\in S_{K_1-1}}\sum_{\frac{a_jmn}{\langle b_jr\rangle_{a_j}+a_ji}<
p\le \frac{a_jln}{\langle b_jr\rangle_{a_j}+a_ji}\atop p\equiv r'\pmod q}\log p+\sum_{p\le
\frac{a_2ln}{\langle b_2r\rangle_{a_2}+a_2K_2}\atop p\equiv r'\pmod q}\log p\\
&\ +\sum_{ \frac{a_1mn}{\langle b_1r\rangle_{a_1}+a_1K_1}<p\le
\frac{a_1ln}{\langle b_1r\rangle_{a_1}+a_1K_1}\atop p\equiv r'\pmod q}\log p+\sum_{i\in S_{K_2-1}\setminus S_{K_1-1}}\sum_{
\frac{a_1mn}{\langle b_1r\rangle_{a_1}+a_1(i+1)}<p\le \frac{a_2ln}{\langle b_2r\rangle_{a_2}+a_2i}
\atop p\equiv r'\pmod q}\log p.
\end{align*}
It then follows from the prime number theorem for arithmetic progressions that

\begin{align*}
&\sum_{p\in \mathcal{P}_{r'}}\log p=\frac{n}{\varphi(q)}
\bigg( \sum_{j=1}^2\sum_{i\in S_{K_1-1}}\Big( \frac{a_jl}{\langle b_jr\rangle_{a_j}+a_ji}-
\frac{a_jm}{\langle b_jr\rangle_{a_j}+a_ji} \Big)+ \frac{a_2l}{\langle b_2r\rangle_{a_2}+a_2K_2}\\
&\quad\quad+\frac{a_1(l-m)}{\langle b_1r\rangle_{a_1}+a_1K_1}+\sum_{i\in S_{K_2-1}\setminus S_{K_1-1}}
\Big( \frac{a_2l}{\langle b_2r\rangle_{a_2}+a_2i}-
\frac{a_1m}{\langle b_1r\rangle_{a_1}+a_1(i+1)}\Big)\bigg)+o(n)\\
&=\frac{n}{\varphi(q)} \bigg(\sum_{i\in S_{K_1}}\frac{a_1l}{\langle b_1r\rangle_{a_1}+a_1i}-
\sum_{i\in S_{K_1-1}}\frac{a_2m}{\langle b_2r\rangle_{a_2}+a_2i}\\
&\quad\quad\quad\quad+\sum_{i\in S_{K_2}}
\Big(\frac{a_2l}{\langle b_2r\rangle_{a_2}+a_2i}-
\frac{a_1m}{\langle b_1r\rangle_{a_1}+a_1i} \Big)\bigg)+o(n)
\end{align*}
as required. Thus (\ref{eq: 4.15}) is true for Case 2.

This completes the proof of Lemma 2.2.
\end{proof}

\section{\bf Proof of Theorem 1.1}

In this section, we use the results presented in Section 2 to give
the proof of Theorem 1.1.

{\it Proof of Theorem 1.1.} Let $r'$ and $r$ be any given integers such that $1\le r', r\le q$ and
$rr'\equiv1\pmod{q}$.  Exchanging $a_1$ with $a_2$ and $b_1$ with $b_2$ simultaneously,
$f(x)=(a_2x+b_2)(a_1x+b_1)$ is unchanged, meanwhile the condition
$a_1\langle b_2r\rangle_{a_2}\ge a_2\langle b_1r\rangle_{a_1}$ in Lemma 2.2 becomes
$a_2\langle b_1r\rangle_{a_1}\ge a_1\langle b_2r\rangle_{a_2}$, and in the conclusion of Lemma 2.2,
$\lambda_r(a_1,a_2,b_1,b_2)$ becomes $\lambda_r(a_2,a_1,b_2,b_1)$.
Thus by Lemma 2.2, we obtain that
$$
\sum_{p\in \mathcal{P}_{r'}}\log p=\frac{n}{\varphi(q)} \lambda_r(a_2,a_1,b_2,b_1)+o(n)
$$
if $a_2\langle b_1r\rangle_{a_1}\ge a_1\langle b_2r\rangle_{a_2}$.
Note that if $a_1\langle b_2r\rangle_{a_2}=a_2\langle b_1r\rangle_{a_1}$, then
$$\frac{a_1mn}{\langle b_1r\rangle_{a_1}+a_1i}=
\frac{a_2mn}{\langle b_2r\rangle_{a_2}+a_2i}$$
and
$$\frac{a_1ln}{\langle b_1r\rangle_{a_1}+a_1i}= \frac{a_2ln}{\langle
b_2r\rangle_{a_2}+a_2i}$$
for any integer $i\ge 0$. Moreover, one has by (\ref{eq: 1.2})
and (\ref{eq: 1.3}) that
$$g_r(a_1, a_2, b_1, b_2)=g_r(a_2, a_1, b_2, b_1)$$
and
$$
h_r(a_1, a_2, b_1, b_2)=h_r(a_2, a_1, b_2, b_1)
$$
if $a_1\langle b_2r\rangle_{a_2}=a_2\langle b_1r\rangle_{a_1}$.
It then follows from (\ref{eq: 1.4}) that
$$\lambda_r(a_1,a_2,b_1, b_2)=\lambda_r(a_2, a_1, b_2, b_1)$$
if $a_1\langle b_2r\rangle_{a_2}=a_2\langle b_1r\rangle_{a_1}$.

Now by Lemma 2.2 and the above discussion, we get that
\begin{align}\label{eq: 4.12}
\sum_{p\in \mathcal{P}_{r'}}\log p=\frac{n}{\varphi(q)}A_r+o(n),
\end{align}
where $A_r$ is defined as in (\ref{eq: 1.5}). Since
$rr'\equiv 1\pmod q$ and $1\le r', r\le q$, $r$ runs over the set
of all positive integers no more than $q$ that are relatively
prime to $q$ as $r'$ does, it then follows from (\ref{eq: 4.12})
and Lemma 2.1 that
\begin{align*}
&\log {\rm lcm}_{mn<i\le ln}\{ f(i)\}\\
=&\frac{n}{\varphi(q)}
\sum_{r'=1\atop \gcd(r',q)=1}^{q}A_r+o(n)\\
=&\frac{n}{\varphi(q)}\sum_{r=1\atop \gcd(r,q)=1}^{q}A_r+o(n)
\end{align*}
as desired. This finishes the proof of Theorem 1.1.
\hfill$\Box$\\

\end{document}